\documentclass{llncs}
\usepackage[left=2.5cm,right=2.5cm,top=2.5cm,bottom=2.5cm]{geometry}
\usepackage{amsmath}
\usepackage{amsfonts}
\usepackage{nicefrac}
\newcommand{\maxent}{q}
\newcommand{\distrFamily}{\mathcal{G}}

\newcommand{\Ex}[1]{\ensuremath{\textstyle E\left(#1\right)}}

\begin{document}
\frontmatter         
\pagestyle{headings}

\title{Maximum Entropy Reconstruction for Discrete Distributions
with Unbounded Support}

\author{Alexander Andreychenko \and Linar Mikeev \and Verena Wolf}

\institute{Saarland University, Saarbrucken, Germany\\
\email{andreychenko@cs.uni-saarland.de}}

\maketitle              % typeset the title of the contribution

\begin{abstract}
The classical problem of moments is addressed by 
the maximum entropy approach for one-dimensional discrete distributions.
The numerical technique of adaptive support approximation is proposed
to reconstruct the distributions in the region where the main
part of probability mass is located.
\end{abstract}

%\category{J.3}{Life and medical sciences}{Biology and genetics}

%\terms{Algorithms, Performance, Theory}

%\keywords{Chemical master equation, moment closure,

%\begin{bottomstuff}
% This research has been partially funded by  the German Research
%     Council (DFG) as part of the Cluster of Excellence
%     on Multimodal Computing and Interaction at Saarland University and
%     the Transregional Collaborative
%     Research Center ``Automatic Verification and Analysis of Complex
%     Systems'' (SFB/TR~14~AVACS).
%
%Author's address: Computer Science Department,
%Saarland University, 66123 Saarbruecken, Germany;
%\end{bottomstuff}

%\maketitle

\section{Introduction}
\label{sec:intro}
In the stochastic chemical kinetics prior information regarding the properties of the 
distribution (e.g. approximately normally distributed) is not directly accessible
and in such a case regaining the probability distribution from the moment description 
is non-trivial. 
In fact 
it turns out that this problem, known as the classical moment problem,
has a long history in other application domains and only 
recently very efficient methods for the reconstruction of the distribution became available.   

Given a number of moments of a random variable, there is in general no 
unique solution for the corresponding distribution.
However it is possible to define
a sequence of distributions that converges to the true one
whenever the number of constraints approaches infinity
~\cite{mnatsakanov_recovery_2009}.  
Conditions for the existence of a solution  
are well-elaborated (such as Krein's and Carlemann's conditions)
but they do not provide a direct algorithmic way to create the reconstruction.
Therefore, Pade approximation
~\cite{mead1984maximum}
and
inverse Laplace transform
~\cite{Chauveau1994186}
have been considered
but turned out to work only in restricted cases and require a
large number of constraints.
Similar difficulties are encountered when 
 lower and upper bounds for the probability distribution are derived
 ~\cite{gavriliadis_moment_2008,tari_unified_2005,Kaas198687}.
Kernel-based approximation methods have been proposed where  one restricts  to a particular class of distributions
~\cite{gavriliadis_truncated_2012,mnatsakanov_recovery_2009,chen_song_2000}.
The numerically most stable methods are, however, based on the 
maximum entropy principle which has its roots 
in statistical mechanics and information theory. The idea is to
choose from all distributions that 
%fulfil
fulfill
the moment constraints the 
distribution that 
%maximises 
maximizes the entropy. 
%Recently there was a large amount of attempts to overcome
%the numerical difficulties arising in the solution process
%and also to incorporate the ability to reconstruct
%multi-dimensional distributions directly.
The maximum entropy reconstruction is the least biased estimate 
that 
%fulfils
fulfills the moment constraints
and it makes no assumptions about the missing information.
No additional knowledge about the
shape of the distribution neither a large number
of moments is necessary.
For instance, if only the first moment (mean) is provided 
the result of applying the maximum entropy principle is 
exponential distribution. In case of two moments
(mean and variance) 
%are provided, 
the reconstruction is given by normal distribution.
%\red{give examples of special cases: first 2 moments gives normal distribution? One moment uniform?}
Additionally, if   experimental data (or simulation traces) is available,
  data-driven 
maximum entropy methods can be applied
~\cite{wu_weighted_2009,Golan1996559}.
Recently, notable progress has been made in the
development of numerical methods for the moment constrained maximum entropy problem
~\cite{abramov2010multidimensional,bandyopadhyay2005maximum,mead1984maximum},
where the main effort is put to the transformation
of the problem in order to overcome the numerical difficulties
that arise during the optimization procedure.
%\red{ one sentence about what they do: rely on a transformation 
%of the optimization problem such that it becomes numerically
%easier to solve?}

In this paper we propose a combination of the classical Newton-based
technique to numerically solve the maximum entropy problem 
with the procedure of distribution support approximation.

\section{Maximum Entropy Reconstruction}\label{sec:maxent}
	The moment closure is usually used to approximate the moments
	of a stochastic dynamical system over time.
	The numerical integration of the correspondent ODE system is
	usually faster than a direct integration of the probability distribution
	or an estimation of the moments based on  Monte-Carlo simulations of the system.
	However, if one is interested in certain events and only the moments of the distribution are known,
	the corresponding probabilities are not directly accessible and
	have to be reconstructed based on the moments.
	Here, we shortly review standard approaches to reconstruct 
	%the 
	one-dimensional 
%	\red{do we restrict to one-dim? We have to say why and how 
%	much effort 2 or more dimensions would be in comparison} 
	marginal probability distributions
	$\pi_i(x_i,t) = P(X_i(t)=x_i)$ of a Markov chain that describes 
	the dynamics of chemical reactions network.
	%according to the theory of stochastic chemical kinetics (cf. Section~\ref{sec:stoch}).
	The task of approximating multi-dimensional distributions 
	follows the	same line however for our case these techniques 
	revealed to be not effective due to numerical difficulties
	in the optimization procedure.
	Thus, we have given (an approximation of) the moments of the $i$-th population and obviously,
	the corresponding distribution is in general not uniquely determined for a finite set of moments.
	In order to select one distribution from this set, we apply the maximum entropy principle.
	~\cite{mead1984maximum}.
	In this way we minimize the amount of prior information about the distribution 
	and avoid any other latent assumption about the distribution.
	Taking its roots in statistical mechanics and thermodynamics
	~\cite{PhysRev.106.620}, 
	the maximum entropy approach was successfully applied 
	to solve moment problems in the	field of 
	climate prediction
	~\cite{abramov2005information,kleeman2002measuring,roulston2002evaluating},
	econometrics
	~\cite{wu_calculation_2003},
	performance analysis
	~\cite{tari_unified_2005,guiasu1986maximum}
	and many others.
%	However, only a few attempts were done to apply it for the 
%	recovery of Markov chains distribution.печковская лиза
%	In this paper we extend the field of applicability to 
%	the systems of biochemical reactions.

	\subsection{Maximum Entropy Approach}
	The maximum entropy principle says that
	among the set of allowed discrete probability distributions $\distrFamily$
 	we choose the probability distribution $q$ that maximizes the entropy $H(g)$ over all  distributions $g \in \distrFamily$, i.e.,
	\begin{equation}
	\label{eq:maxShannonProblem}
	\begin{array}{c}
		\maxent = \arg \max_{g\in \distrFamily} H(g)
				= \arg \max_{g\in \distrFamily} \left( -\sum_x g(x) \ln{g(x)}\right)
				.
	\end{array}
	\end{equation}
	where $x$ ranges over all possible states of the discrete state space.
	Note that we assume that all distributions are defined on the same state space.	In our case the set $\distrFamily$
	consists of all discrete probability distributions that satisfy the moment constraints. Given a sequence of 	$M$ non-central moments
	%\red{index $i$ is used for enumerating species! Better use $k$ or so}
	%$$\Ex{X^i}=\mu_i, i=0,1,\ldots,m,$$
	$$\Ex{X^k}=\mu_k, k=0,1,\ldots,M,$$
	 the following constraints are considered
	%$$\int\limits_{D} x^i \maxent(x) dx = \mu_i$$
	\begin{equation}\label{eq:momentconstr}
	\sum_x x^k g(x)   = \mu_k, k=0,1,\ldots,M.
	\end{equation}
	%\fxnote{where $D$ is the support of $f(x)$.}.
%	Here, we add the constraint for $\mu_0=1$ in order to choose from
%	a set of non-negative functions and 	ensure that $g$ is a distribution.
	Here, we choose $g$ to be a non-negative function 
	and add the constraint $\mu_0=1$ in order 
	to ensure that $g$ is a distribution.
	The above problem is a nonlinear constrained optimization problem, which is usually
	addressed by the method of Lagrange. Consider the Lagrangian
	functional
	\begin{equation*}
	\begin{array}{c}
		\mathcal{L}(g,\lambda)
		= H(g) - \sum\limits_{k=0}^{M} \lambda_k
				\left( \sum_x x^k g(x)  - \mu_k \right),
	\end{array}
	\end{equation*}
	where $\lambda=(\lambda_0,\ldots,\lambda_M)$ are the corresponding Lagrangian multipliers.
	It is possible to show that maximizing the  unconstraint Lagrangian $\mathcal{L}$
	gives a solution to the constrained maximum entropy problem
	%~\ref{eq:maxShannonProblem}.
	%subject to the constraints in Eq.~\eqref{eq:momentconstr}
 	The variation of the functional $\mathcal{L}$ according to
	the unknown distribution provides the general form of $g(x)$
	$$	\frac{\partial \mathcal{L}}{\partial g(x)} = 0
		\implies
		g(x) = \exp \left( -1 -\sum\limits_{k=0}^{M} \lambda_k x^k \right)
		=\frac{1}{Z(x)} \exp \left( -\sum\limits_{k=1}^{M} \lambda_k x^k \right),
	$$
	where
	%$Z(x) = e^{\lambda_0} = \int \exp \left( \sum\limits_{i=0}^{m+1} \lambda_i \mu_i \right)$
	%$$Z(x) = e^{-1-\lambda_0}
	%\change{
	\begin{equation}
		\label{eq:normalization_constant_Z}
			Z(x) = e^{1+\lambda_0}
			      = \displaystyle\sum_x \exp \left( -\sum\limits_{k=1}^{M} \lambda_k x^k \right)
	\end{equation}
%	$$
%	\label{eq:normalization_constant_Z}
%	Z(x) = e^{1+\lambda_0}
%	      = \displaystyle\sum_x \exp \left( -\sum\limits_{k=1}^{M} \lambda_k x^k \right) $$
	%      }
	%\todo{reviewers did not notice the mistake here but it has to be fixed}
	is a normalization constant.
	In %order to cope with the constrained optimization problem for $\mathcal{L}$,
	the dual approach
	%\red{~\cite{judge1996maximum} - really needed?}
	we insert the above equation for $g(x)$ into
	the Lagrangian thus we can transform the problem into an
	unconstrained convex minimization problem of the dual function w.r.t
	to the dual variable $\lambda$
	%$$\Psi(\lambda)=\ln Z(x) - \sum\limits_{i=1}^{M} \lambda_i \mu_i,$$
	$$\Psi(\lambda)=\ln Z(x) + \sum\limits_{k=1}^{M} \lambda_k \mu_k,$$
		According to the Kuhn-Tucker theorem, the solution
	$\lambda^* = \arg \min \Psi(\lambda)$ of the minimization problem
	for the dual function equals the solution $q$ of the original constrained optimization problem~\eqref{eq:maxShannonProblem}.

	\subsection{Maximum Entropy Numerical Approximation}
	It is possible to solve  the constrained maximization problem in Eq.~\eqref{eq:maxShannonProblem}
	%for $m \leq 2$ analytically. 
	for $M \leq 2$ analytically. 
	%For $m>2$ 
	For $M>2$ 
	%see ... for the existence conditions)
	numerical methods have to be applied
	to incorporate the knowledge of  moments of order three and more.
	%\change{
	Here we use the Levenberg-Marquardt 
	method
	~\cite{transtrum_improvements_2012} 
	to minimize the dual function 
	$\Psi(\lambda)$.
	%Approximate solution $\hat{q}$ to the problem~\eqref{eq:maxShannonProblem}
	An approximate solution $\tilde{q}$ is given by
	%$$ \hat{q}(x) = \exp \left( -1 -\sum\limits_{k=0}^{M} \hat{\lambda}_k x^k \right) , $$
	$$ \tilde{q}(x) = \exp \left( -1 -\sum\limits_{k=0}^{M} \hat{\lambda}_k x^k \right) , $$
	%where $\hat{\lambda}$ is the result of iterative procedure ~\eqref{eq:dualFunctionNewtonOpt}
	where $\tilde{\lambda}$ is the result of the iteration 
	\begin{equation}
	\label{eq:dualFunctionNewtonOpt}
			\lambda^{(\ell+1)} = \lambda^{(\ell)} - \left( H + \gamma^{(\ell)}
		 \cdot \mathrm{diag}(H) \right)^{-1} \frac{\partial \Psi}{\partial \lambda}.
		\end{equation}
	The damping factor $\gamma$ is updated
	according to the strategy suggested in
	~\cite{transtrum_improvements_2012}
	and $\lambda^{(\ell)}=(\lambda^{(\ell)}_1,\ldots,\lambda^{(\ell)}_M)$
	is an approximation of 
	the vector $\lambda=(\lambda_1,\ldots,\lambda_M)$ in the $\ell$-th step of the
	iteration. We compute $\lambda_0$  
	as $\lambda_0 = \ln Z - 1$ (see Eq.~\eqref{eq:normalization_constant_Z}).	
	Initially we choose
	$\lambda^{(0)}=(0,\ldots,0)$ and stop when
	the solution converges, i.e. when the condition
	$ \vert \lambda^{(\ell+1)} - \lambda^{(\ell)} \vert < \delta_\lambda $
	is satisfied for a small threshold $\delta_\lambda>0$.
	In the $\ell$-th iteration  the components of the gradient vector 
	are approximated by
	$\frac{\partial \Psi}{\partial \lambda_i} %= 
			\textstyle\approx \mu_i - \frac{1}{Z} \widetilde{\mu}_i$
	and the entries of the Hessian matrix are computed as
	$$
		H_{i,j} = \frac{\partial^2 \Psi}{\partial \lambda_i \partial \lambda_j} 
							\approx \frac{Z \cdot \widetilde{\mu}_{i+j} - \widetilde{\mu}_i \widetilde{\mu}_j}{Z^2},
							\quad i,j=1, \ldots, M.
	$$
	The approximation $\widetilde{\mu}_i$ of the $i$-th moment is given by 
		\begin{equation}
		\label{eq:momentapp}
		\textstyle	\widetilde{\mu}_i 
			= 
			\sum\limits_{x} 
			x^i \exp \left( -\sum\nolimits_{k=1}^{M} \lambda_k^{(\ell)} x^k \right),
			\quad i=1,\ldots,2M,
		\end{equation}
	%Initially we assume that all discrete distributions $g \in \distrFamily$
	%are defined on $\mathbb{Z}^+$.
	In order to approximate the moments 
	we need to truncate the infinite sum in  Eq.~\eqref{eq:momentapp}.
	We refer to Section~\ref{app:B} for a detailed description 
	   of how the distribution support can be
			  approximated.
			
	The convexity~
	\cite{mead1984maximum} 
	of the dual function $\Psi(\lambda)$ 
	guarantees the existence of a unique minimum $\lambda^*$ 
	%approximated by $\hat{\lambda}$.
	approximated by $\tilde{\lambda}$.
	Theoretical conditions for the existence of the solution 
	are discussed in detail in~
	\cite{tari_unified_2005,Stoyanov_2000,Lin199785}.
	A similar analysis for the multivariate case is provided
	in~ 
	\cite{kleiber_multivariate_2013}.
	The iterative procedure in Eq.~\eqref{eq:dualFunctionNewtonOpt}
	might however fail due to   numerical instabilities
	when the inverse of the Hessian is calculated.
	The iterative minimization presented in~
	\cite{bandyopadhyay2005maximum}
	and the Broyden-Fletcher-Goldfarb-Shanno (BFGS) 
	procedure~
	\cite{byrd1995limited}
	can be used to improve the numerical stability.	
	%In the sequel we denote by $\hat{q}_i(x)$
	In the sequel we denote by $\tilde{\pi}_i(x,t)$
	the reconstructed distribution of the $i$-th species
	for a given sequence of moments $\mu_0, \ldots, \mu_M$,
	i.e. the marginal probability distribution 
	%$\pi_i(\cdot,t) = P \left( X_i = x \right)$
	$\pi_i(x,t) = P \left( X_i = x \right)$.
	Note that the reconstruction approach presented above
	provides a reasonable approximation of the probabilities
	only in high-probability regions.
		In order to accurately approximate the tails of the distribution 
	special	methods have been developed~\cite{gavriliadis_truncated_2012}.
	%}

\section{Conclusions}\label{sec:conc}
As future work, we plan to extend the reconstruction procedure 
in several ways. First, we want to consider moments of higher order than five.
Since in this case the concrete values become very large it might be advantageous 
to consider central moments instead which implies that the reconstruction procedure
has to be adapted.
Alternatively, we might (instead of algebraic moments) consider 
other functions of the random variables such as
exponential functions~\cite{mnatsakanov_note_2013},
Fup functions~
\cite{gotovac_maximum_2009},
and
Chebyshev polynomials~
\cite{bandyopadhyay2005maximum}.
 Another possible extension could
address  the problem of truncating the support of the distribution
such that the reconstruction is applied to a finite support. 
We expect that in this case the reconstruction will become more
accurate since we will not have to rely on the Gauss-Hermite quadrature formula.
%
%direct reconstruction
%of the discrete distribution with semi-infinite support,
%where not only the maximum entropy principle can be applied
%but other methods of solving the classical moments as well
%(cf. Section~\ref{sec:intro}).
For instance, the theory of Christoffel functions~
\cite{gavriliadis_truncated_2012}
could be used to determine the region where
the main part of the probability mass is located. 

Finally, we want to improve the approximation for species that are present
in very small quantities, since for those species a direct representation of the probabilities 
is more appropriate than a moment representation. Therefore  we plan to consider 
the conditional moments approach~\cite{MCM_Hasenauer_Wolf}, 
where we only integrate the moments of species having large molecular counts but 
keep the discrete probabilities   for the species with small populations.

% Acknowledgments
%\begin{acks}
%This research has been partially funded by  the German Research
%     Council (DFG) as part of the Cluster of Excellence
%     on Multimodal Computing and Interaction at Saarland University and
%     the Transregional Collaborative
%     Research Center ``Automatic Verification and Analysis of Complex
%     Systems'' (SFB/TR~14~AVACS).
%     \end{acks}

%\newpage
\begin{appendix}

\section{Approximation of the Support}
\label{app:B}
\newcommand{\maxEntProbThreshold}{\delta_{\mbox{\scriptsize prob}}}
	%\change{
	During the iteration procedure~\eqref{eq:dualFunctionNewtonOpt}
	we need to approximate one-dimensional moments by summing up over all
	states $x \in \mathbb{Z}^+$ that have positive probability mass.
	However our case studies possess the infinite number of such states
	and the appropriate truncation has to be done.
	Instead of considering whole state space $\mathbb{Z}^+$
	we consider a subset $D = \{x_L,\ldots,x_R\}\subset \mathbb{Z}^+$,
	where we have to choose such values for $x_L$ and $x_R$
	%for which the iteration procedure converges. 
	that the iteration procedure converges. 
	It might fail to converge if the difference $(x_R - x_L) $
	is very large so that the conditional number of the matrix
	$ \left( H + \gamma^{(\ell)} \cdot \mathrm{diag}(H) \right) $
	is very large.
	To find a reasonable initial guess 
	$D^{(0)} = \{ x_L^{(0)}, \ldots,  x_R^{(0)} \}$
	we use the results in~\cite{tari_unified_2005}
	% to find a region
	%	that contains the main part of the probability mass.
	and consider the roots of the function $\Delta^0_k(w)$
	\begin{equation}
	%\label{eq:delta0}
		\Delta^0_k(w) = 
				\left|
					\begin{array}{cccc}
						\mu_0 & \mu_1 & \cdots & \mu_k \\
						\vdots & &  & \vdots \\
						\mu_{k-1} & \mu_k & \cdots & \mu_{2k-1} \\
						1 & w & \cdots & w^k \\
					\end{array}
				\right|,
	\end{equation}
	where $k=\lfloor \frac{M}{2} \rfloor$, and $M$ is even.
	The initial guess $D^{(0)}$ is defined by 
	$x_L^{(0)} = \lfloor w_1 \rfloor$ and $x_R^{(0)} = \lceil w_k \rceil$,
	where $w_1 < \ldots < w_k$ are real and simple roots
	of the equation $\Delta^0_k(w) = 0$.
	In the $\ell$-th iteration we check 
	if the probability of the right-most state $x_R^{(\ell)}$
	is reasonably small in comparison to the 
	maximum value of $\tilde{q}^{(\ell)}(x)$ for $x\in D^{(\ell)}$, i.e.  
	\begin{equation}
	\label{eq:supportInequality}
			\tilde{q}^{(\ell)}(x_R^{(\ell)}) 
			< \maxEntProbThreshold \cdot \max\nolimits_{x \in D^{(\ell)}} \tilde{q}^{(\ell)}(x),
	\end{equation}
	where $\maxEntProbThreshold$ is a small threshold
	(for all out experiments we chose $\maxEntProbThreshold = 10^{-3}$).
	We extend the support until inequality~\eqref{eq:supportInequality}
	is satisfied by adding new state on each iteration
	\begin{equation}
		\left( x_L^{(\ell+1)}, x_R^{(\ell+1)} \right) = 
		\begin{cases}
		  \left( \max (0, x_L^{(\ell)}), x_R^{(\ell)} \right), & \ell \text{  is even} \\
   	  \left( x_L^{(\ell)}, x_R^{(\ell)}+1 \right), & \ell \text{  is odd}.
		  \end{cases}
	\end{equation}
	The final results $\tilde{\lambda}$ and $\hat{D}$ of the iteration yields the distribution 
	$\tilde{q}(x)$ that approximates the marginal distribution of interest.

	Please note that $M$ is assumed to be even when we use the function $\Delta_k^0$.
	Tari et. al also provide the extension of this technique that allows to
	account for the case when odd number of moments is known ($M$ is odd)
	%In addition to the function $\Delta^0(w)$ defined in Eq.~\eqref{eq:delta0},
	by considering the function $\Delta^1_z(\eta)$
	\begin{equation*}
			\Delta^1(\eta) = 
				\left|
					\begin{array}{cccc}
							\mu_1 - w_1 \mu_0 & \mu_2-w_1 \mu_1 &\cdots & \mu_{z} - w_1 \mu_{z-1} \\
							\vdots & \vdots &  & \vdots \\
							\mu_{z-1} - w_1 \mu_{z-2} & \mu_{z} -w_1 \mu_{z-1} & \cdots & \mu_{2z-2} \! -\! w_1 \mu_{2z-3} \\
							1 & \eta & \cdots & \eta^{z-1} \\
					\end{array}
				\right|,
	\end{equation*}
	where $z=\lfloor \frac{M}{2} \rfloor + 1$.
	%and $w_1$ is 
	%the smallest root of the equation $\Delta^0(w)=0$.
	Let $W = \{w_1, \ldots, w_k \}$ be the set of the solutions
	of $\Delta^0_k(w)=0$ 
	and 
	$H= \{ \eta_1, \ldots, \eta_z \}$ be
	the set of solutions of $\Delta^1_z(\eta) = 0$,
	where all the elements of $W$ and $H$ are real and simple.
	The first approximation $D^{(0)}$	is then defined by
	$x_L^{(0)} = \lfloor \min(w_1, \eta_1) \rfloor$
	and $x_R^{(0)} = \lceil \max(w_k, \eta_z) \rceil$.
	
	\paragraph{Alternative method for support approximation extension.}
	Instead of adding only one state to the support approximation $D^{(\ell+1)}$
	we can use the following heuristics on 
		Chebyshev's inequality 
		$P\{ \vert X - \mu_1 \vert \geq z \} \leq \xi$,
		$\xi = \nicefrac{\mu_2}{z}$.
		%as follows.
		We compute $z=\nicefrac{\mu_2}{0.1}$ (where we fix $\xi = 0.1$) 
		such that the correspondent set 
		$\widetilde{D} = \{ \lfloor \mu_1-z \rfloor, \ldots, \lceil \mu_1 + z \rceil \}$
		shall include at least $90\%$ of the probability mass.
		\newline
		The initial approximation 
		%for the support 
		$D^{(0)}$ is compared to $\widetilde{D}$ 
		by computing the difference $l = \vert z - x_R^{(0)} \vert$.
		The latter serves as an increment in support extension procedure,
		i.e. the approximation of the support on iteration $(i+1)$
		is given by $D^{(i+1)} = \{ x_L^{(i+1)}, \ldots,  x_R^{(i+1)} \}$,
		where
		$x_L^{(i+1)} = \max \left( 0, x_L^{(i)} - \lceil \nicefrac{l}{2} \rceil \right)$
		and 
		$x_R^{(i+1)} = x_R^{(i)} + \lceil \nicefrac{l}{2} \rceil$.
	
\end{appendix}

% Bibliography
\bibliographystyle{splncs03}
\bibliography{techrep}
%\clearpage 

\end{document}